\documentclass[12pt]{amsart}

 \usepackage{amsmath}
\usepackage{amssymb,latexsym,cite}
\usepackage{amscd}
\usepackage{comment}
 \usepackage{xspace}
 \usepackage{verbatim}
 \usepackage{todonotes}

\newtheorem{theorem}{Theorem}[section]

\newtheorem{lemma}[theorem]{Lemma}
\theoremstyle{definition}
\newtheorem{remark}[theorem]{Remark}

\numberwithin{equation}{section}

\newcommand{\R}{{\mathbb R}}

\newcommand{\BSA}{\begin{subarray}}
\newcommand{\ESA}{\end{subarray}}
\newcommand{\BAL}{\begin{aligned}}
\newcommand{\EAL}{\end{aligned}}

\newcommand{\forevery}{\quad \forall}

\newcommand{\rec}[1]{\frac{1}{#1}}

\newcommand{\dist}{\mathrm{dist}\,}

\newcommand{\prt}{\partial}
\newcommand{\sms}{\setminus}

\newcommand{\ti}{\times}

\newcommand{\tl}{\tilde}

\newcommand{\sth}{such that\xspace}

\newcommand{\bvp}{boundary value problem\xspace}

\newcommand{\bdw}{\partial\Gw}

\newcommand{\qtxt}[1]{\quad\textrm{#1}}

     \def\gb{\beta}       \def\gg{\gamma}
       \def\gd{\delta}      \def\ge{\epsilon}

      \def\gk{\kappa}      
                 
    \def\gr{\rho}        
       
      \def\gw{\omega}
\def\gx{\xi}                
\def\Gg{\Gamma}     \def\Gd{\Delta}

\def\Gw{\Omega}              

      \def\CF{{\mathcal F}}

   \def\BBR {\mathbb R}

\def\rs{radially symmetric\xspace}

\begin{document}

\title[Large solutions]{Large solutions of semilinear equations with Hardy potential}
\author{Moshe Marcus}
\address{Department of Mathematics\\ Technion\\ Haifa 32000\\ ISRAEL}
\email{marcusm@math.technion.ac.il}

\date{\today}

\maketitle

\section{Introduction} Assume that $\Gw$ is a bounded Lipschitz  domain in $\R^N$. Denote
\begin{equation}\label{Lmu}
L_\mu=\Gd+\mu/\gd^2, \quad \gd(x)=\dist(x,\bdw),\quad \mu\in \R.
\end{equation}
\[ -L_\mu u+f(u)=0, \tag{NE}\]
\[  -\Gd u+f(u)=0. \tag{NE0}\]
Assumptions on $f$:
\begin{equation}\label{con0}
\BAL f\in C^1[0,1), \;& f(0)=0, \quad f \text{ and } f'\;   \text{positive on $(0,\infty)$}\\
& h(t):=f(t)/t\to \infty \qtxt {sa }\; t\to\infty.\EAL
  \tag{F1}
\end{equation}

\noindent(F2)\hskip 4mm  $\exists a^\#>0$ \sth, for every positive solution $u$ of (NE0),
\begin{equation}\label{con1}
h(u(x))\leq a^\#\gd(x)^{-2} \forevery x\in \Gw,
\end{equation}
where $h(t):=f(t)/t, \forevery t>0$.

Note that if $f$ satisfies (F1) and is convex on $\R_+$ then $h$ is non-decreasing.

The Keller -- Osserman condition on $f$ is a key condition in the study of equations such as (NE0). The condition is:
\[ \psi(t):=\int_t^\infty (2F(s)^{-1/2}<\infty \forevery t>0 \tag{KO} \]
where $F(s)=\int_0^s f$. Since  (F2) implies that the set of positive solutions of (NE0) is locally uniformly bounded, conditions (F1) and (F2) imply the Keller -- Osserman condition (se \cite{K} and \cite{BM92}).The converse does not hold in general, but it is still valid under an additional condition on $f$ \cite[Lemma 5.1]{MM-nl}:

\begin{lemma}\label{VR}
 Assume that $f$ satisfies (F1) and the Keller - Osserman condition (KO). Let $V_R$ denote the large solution of (NE0) in $B_R(0)$.
Given a number $c_0>0$, suppose that there exists a constant $C_1>0$ \sth
 \begin{equation}\label{C1}
 \psi(2s) \leq C_1h(s)^{-1/2} \qtxt{whenever }\;  h(s)\geq c_0.
 \end{equation}
 Then
 \begin{equation}\label{C2}
 h(V_R(0))\leq C_2R^{-2} \qtxt{where }\;  C_2= N(2+C_1)^2.
  \end{equation}
\end{lemma}

 The  existence of $c_0>0$ \sth \eqref{C1} holds depends  only on $f$. Consequently $C_2$ depends only on $f$, $N$ and the choice of $c_0$.

Condition \eqref{C1} holds for a large family of functions $f$ including $f(t)= t^p$, $p>1$ and $f(t)=e^t-1$.


More general results dealing with (NE0) when $f$ also depends  on the space variable were obtained by Ancona in an appendix to \cite{AM}. Naturally in that case, in general, the constant $C_2$ depends also on the domain.

 Further, by \cite[Theorem 6.1]{MM-Tai2}, (F2) implies that there exist positive constants $a_0$, $a_1$ such that, for every positive solution u of (NE) in $\Gw$,
\begin{equation}\label{F2'}
 h(u(x)/a_1) \leq a_0\gd(x)^{-2} \forevery x\in \Gw. \tag {F2'}
\end{equation}
The constants $a_0$, $a_1$  depend on the constant in \eqref{C2} and the constant in the strong Harnack inequality associated with $L_\mu$.


If $\Gw$ is a bounded $C^2$ domain, $f$ is convex and  satisfies conditions (F1) and (KO) then there exists a unique large solution of (NE0), say $U^\Gw_f$, and it satisfies,
\begin{equation}\label{u/phi}
\frac{U^\Gw_f(x)}{\phi(\gd(x))}\to 1 \qtxt{as }x\to\bdw,
\end{equation}
where $\phi$ is the inverse of the  function  $\psi$ defined in (KO).  This result was first proved in \cite{BM92} under a stronger assumptions on $f$. Under the conditions stated above, it  was established in \cite{BM98}
(see Theorem 3 and inequality (34) there).


The function $\phi$ is  the unique solution of the problem,
\begin{equation}\label{1d}
\phi''= F(\phi) \qtxt{in }\R_+, \quad \lim_{t\to 0}\phi(t)=\infty.
\end{equation}

In fact, under the same assumptions on $f$, uniqueness of the large solution  (but not \eqref{u/phi}) holds in much more general domains. In particular,
if $\Gw$ is a bounded Lipschitz domain, $f$ is convex and  satisfies  (F1) and (KO) then there exists a unique large  solution of (NE0) \cite[Thm. 1.4]{MV2006}.

In this note we discuss the existence and uniqueness of large solutions of equation (NE) when $\mu\geq 0$ under various conditions on $f$ and on the domain. We establish the following results.

Denote by $C_H(\Gw)$ the Hardy constant in $\Gw$ relative to the potential $\gd^{-2}$ i.e.
$$ c_H(\Gw)=\inf_{\varphi\in C_c^\infty(\Gw)}\frac{\int_\Gw|\nabla \varphi|^2}{\int_\Gw \gd^{-2}\varphi^2}.$$

\begin{theorem}\label{exist}
Let $\Gw$ be a bounded Lipschitz domain. Assume that $f$ satisfies (F1) and (F2) and that it is convex on $\R_+$. Then
for every $\mu\geq 0$, there exists a large solution $U$ of (NE) \sth $U>U_f^\Gw$.
\end{theorem}

\begin{theorem}\label{uniq-1}
	Let $\Gw$ be a bounded Lipschitz domain. Assume that $f$ satisfies (F1) and (F2) and that it is convex on $\R_+$.
	
		Under these assumptions, if $0\leq\mu < 1/4$ then (NE) has a unique large solution.
	\end{theorem}

\begin{theorem}\label{Main}
Let $\Gw$ be a bounded $C^2$ domain. Assume that $f$ is convex and satisfies conditions(F1), (F2). In addition assume that:\\ 

 For every $a>1$ there exist $c>1$ and $t_0>0$ \sth
\begin{equation}\label{con2}\BAL
 ah(t)\leq h(ct), \quad t>t_0
\EAL\end{equation}

 For every $a\in (0,1)$ there exist $c\in (0,1)$ and $t_0>0$ \sth
\begin{equation}\label{con2.0}
h(c  t)\leq a h(t), \quad t>t_0.
\end{equation}

There exists $A>1$ \sth
\begin{equation}\label{con3}
 h(\phi)\leq  A\gd^{-2}.
\end{equation}

Then, for every $\mu>0$,  equation (NE) has a unique large solution $\bar U_\mu$. Moreover there 
 exist positive constants $c_1, c_2$ \sth
\begin{equation}\label{U_bounds}
 c_2 \phi \leq \bar  U_\mu \leq c_1\tl\phi \qtxt{where }\; \tl\phi:=h^{-1}(\gd^{-2}).
\end{equation}
\end{theorem}


\remark $\,$ The conditions of Theorem \ref{Main} are satisfied in particular by the functions $f(u)= u^p$, $p>1$ and $f(u)=e^u-1$. 

If $f$ is a power or, more generally, 
\begin{equation}\label{powers}
b_2\leq  f(t)/t^p \leq b_1 \forevery t>t_0 
\end{equation} 
for some positive constants, $b_1, b_2$, $t_0$
then, 
\[\exists C>0: \quad \bar \phi \leq C\phi. \tag{\#}\]

%

Equations of the form (NE) have been studied intensively in the last decade. Among the earliest were the works of Bandle, Moroz and Reichel \cite{BMR1} and \cite{BMR2}. The first of these papers dealt with power nonlinearities. The second studied the case $f(u)=e^u$ and proved existence and uniqueness of the large solution, in smooth domains,   when $0<\mu<c_H(\Gw)$ (= the Hardy constant in $\Gw$). Recall that, in general $c_H(\Gw)\leq 1/4$.

 More recent works include Marcus and P.T. Nguyen \cite{MM-Tai1} and \cite{MM-Tai2} and Gkikas and Veron \cite{Gk+V}. These works  dealt with boundary value problems for (NE) with power nonlinearities, in $C^2$ domains.
Under the same conditions,  equation (NE) was also studied by Du and Wei  \cite{Du_Wei1} who proved that, in this case, (NE) has a unique large solution for every $\mu>0$.

Recently, some  questions about large solutions of (NE) in a ball, with $\mu<1/4$ and $f(u)=u^p$, $p>0$ have been studied by Bandle and  Pozio \cite{Ba_Po}.

Existence and uniqueness of large solutions for the equation
$$-(\Gd u +\mu/|x|^2)u + u^p=0, \quad p>1$$
in smooth domains $\Gw$ (with $0\in \bdw$)
 have been studied by  Guerch and Veron \cite{Gu+V}, F.C. Cirstea \cite{Cir}, Du and Wei \cite{Du_Wei2} a.o.


%


\section{Proof of Theorem \ref{exist}}
Under the assumptions of the theorem, equation (NE0) possesses a unique large solution, say $U_f^\Gw$, (Marcus and Veron \cite[Thm. 1.4]{MV2006}).
Since $\mu>0$, $U_f^\Gw$ is a subsolution of (NE).

Let $\{D_n\}$ be a smooth exhaustion of $\Gw$ and denote by $u_n$ the solution of the \bvp
$$\BAL -L_\mu u+f(u)&=0 \qtxt{in }\; D_n\\
u&=U_f^\Gw  \qtxt{on }\; \prt D_n.
\EAL $$
Then $\{u_n\}$ is a monotone increasing sequence.  For every set $E\Subset \Gw$ let $n_E$ be a number \sth $E\Subset D_n$ for all $n\geq n_E$. Condition (F2')implies that $\{u_n:n\geq n_E\}$  is uniformly bounded in $E$. Therefore the sequence converges to a solution $U$ of (NE).  Since $U>U_f^\Gw$ it follows that $U$ is a large solution of (NE).

\qed

\remark The above argument shows tat if $u$ is any positive solution of (NE0) there exists a solution $\tl u$ of (NE) \sth $\tl u>u$. Moreover, assuming (F1) and (KO), if $\Gw$ is a bounded domain \sth (NE0) has a large solution then (NE) has a large solution.

\section{Proof of Theorem \ref{uniq-1}}

We start with some notation.
Denote,
$$T(r,\gr)= \{\gx= (\gx_1, \gx')\in \R\ti\R^{N-1}:    |\gx_1|<\gr,\;  |\gx'|<r\}.$$
By assumption $\Gw$ is a bounded Lipschitz domain. Consequently there exist positive numbers  $r_0$ , $\gk$ \sth, for every $y\in \bdw$, there exist: (i) a  set of Euclidean coordinates $\gx=\gx_y$ centered at $y$ with the positive $\gx_1$ axis pointing in the direction of the inner normal $\mathbf{n}_y$ and (ii) a Lipschitz function $F_y$ on $\R^{N-1}$ with Lipschitz constant $\leq \kappa$ \sth $F_y(0)=0$,  and
\begin{equation}\label{ngh_y}\BAL
Q^y(r_0,\gr_0):=\, &\Gw\cap T^y(r_0,\gr_0)\\
=\, &\{ \gx= (\gx_1, \gx'): F^y(\gx')<\gx_1<\gr_0,\; |\gx'|<r_0\},
\EAL\end{equation}
where $\gx=\gx^y$, $ T^y(r_0,\gr_0)= y  +  T(r_0,\gr_0)$ and $\gr_0 = 10\gk r_0$.
Without loss of generality, we assume that $\kappa>1$.

The set of coordinates $\gx^y$ is called a standard set of coordinates at $y$ and $T^y(r,\gr)$ with $0<r\leq r_0$ and $\gr=c\gk r$, $2<c\leq 10$ is called a standard cylinder at $y$.

Suppose there exist two large solution $u_1, u_2$. We may assume that $u_1\leq u_2$. Otherwise we replace $u_2$ by the solution lying between the subsolution $\max(u_1,u_2)$ and the supersolution $u_1+ u_2$.
Here we use the fact that, as $f$ is convex, $f(0)=0$ and $f$ is increasing,
$$ f(a) + f(b)\leq f(a+b) \forevery a,\,b >0.$$

In what follows $y$ is kept fixed and we drop the superscript in $Q^y$. Further, for $0<a\leq 1$ we denote
 $$aQ=Q(ar_0,a\rho_0), \quad \Gg_ {1,a}=\prt (aQ)\cap\bdw, \quad
\Gg_{2,a}=\prt (aQ) \cap\Gw.$$


\noindent\textbf{Part 1.}
 \textit{Construction  of a \emph{subsolution} $w$ of (NE) in $\gw:=\rec{2}Q$ \sth:} \\ [2mm]
 	\noindent 	(i) $w\in C(\bar \gw \cap\Gw)$, \\ [1mm]
(ii)	$w=0$ on $\prt \gw \cap\Gw$,\\ [2mm] 
(iii) $\dfrac{w}{u_2}\to 1$ as $\xi\to \Gg_ {1,a}$  uniformly  for $|\xi'|<ar_0$,
  $\forevery a\in (0,\rec{2})$.\\ [1mm]

Let
$$ \gw_n:=\{\xi\in\gw: F(\xi')+\frac{\rho_0}{8n} <\xi_1\}.$$


Let $v_n$ be the solution of the \bvp,

$$\BAL -L_\mu v_n & =0 \qtxt{in }\; \gw_n,\\
v_n & =u_2-j_n \qtxt{on }\; \prt\gw_n,
\EAL $$
where $j_n$ is a non-negative, continuous function  on $\prt\gw_n$ \sth
$$j_{n+1}\leq j_n\leq u_2$$
and
$$ j_n=\begin{cases} u_2 &\qtxt{on }\Gg'_{n}=\prt\gw_n\cap \{\xi_1=F(\xi')+\frac{\rho_0}{8n}\},\\
0 & \prt\gw_n\cap \{\xi_1>F(\xi')+\frac{\rho_0}{4n}\}.
\end{cases}
$$
In particular $v_n=0$ on $\Gg'_{n}$ and $v_n\leq u_2$ in $\gw_n$.

Let $G_{\mu,n}$ be the Green function
of $L_\mu$ in
$$Q_n :=\{\xi\in \frac{3}{4}Q: F(\xi')+\frac{\rho_0}{8n} <\xi_1\}.$$
and let $\gx_0=(\frac{2}{3}\rho_0, 0)$.

By the Boundary Harnack principle applied to $L_\mu$ in $\gw_n$,  for every $a\in (0,1/2)$ there exists a constant $C_a$ \sth,
$$v_n(\xi) \leq C_aG_{\mu,n}(\xi,\xi_0) \qtxt{in }\;\gw_{n,a}:= \{\xi\in aQ: F(\xi')+\frac{\rho_0}{8n} <\xi_1\}.$$


Recall that, if $\Gw$ is Lipschitz and $\gr$ is sufficiently small then, by \cite{MMP},  the local Hardy constant in the strip  $\Gw_\gr:=\{x\in \Gw: \gd(x)<\gr\}$ 
is $1/4$, although the Hardy constant in $\Gw$ may be lower. Therefore, for $\mu<1/4$, $L_\mu$ has a Green function $G_{\mu,Q}$ in $Q$. 
Since $G_{\mu,n}$ is dominated by $G_{\mu,Q}$ 
and $G_{\mu,Q}(\xi, \xi_0)$ is bounded in $\gw$
it follows that there exists a constant $C'(a)$  \sth,
\[\sup_{\gw_{n,a}} v_n\leq C'_a \forevery a\in (0,1/2).\tag{1*}\]
As $\{v_n\}$ increases, the sequence converges to a solution $v$ of $L_\mu v=0$ in $\gw$ and the convergence is uniform in $aQ$ for every $a\in (0,1/2)$.
In addition $\{v_n\}$ is bounded by $u_2$ which is continuous in
$$ \{\xi\in \bar \gw : \xi_1>\ge\} \forevery \ge\in (0,\rho/2).$$
Therefore, for every $\ge$ as above, the sequence converges uniformly in this set.  Thus $v$ is continuous in $\{\xi\in \bar \gw,\; \xi_1>\ge\}$   and
\[ v=u_2 \qtxt{on } \prt\gw\cap \Gw. \tag{2*}\]

Put $ w=u_2 - v$. Then $w$ is a subsolution of (NE) in $\gw$:
\[-L_\mu w + f(w)< -L_\mu u_2 +f(u_2)=0. \tag{3*}\]
 By (1*)
$$v=u_2-w\leq C'(a)\qtxt{in }aQ.$$
By assumption, $u_2\to\infty$ as $\xi\to \prt \gw\cap\bdw$. Therefore
\[  1-\frac{w}{u_2}\to 0 \qtxt{as }\;\xi\to \Gg_{1,a}, \forevery a\in (0,1/2). \tag{4*}\]

In conclusion $w$ has the properties stated in Part 1.\\ [1mm]


%

\noindent\textbf{Part 2.} \textit{Completion of proof.}\\ [2mm]
 Given $\gb\in (0,\rho_0/10)$  denote
$$w^\gb(\xi)=w(\xi_1+\gb,\xi')  \qtxt{in }\; \gw^\gb:=\{\xi\in \gw:\,\xi_1<\frac{\rho_0}{2}-\gb\}.$$
By (3*), $w^\gb$ is a subsolution,
$$ -(\Gd+ \frac{\mu}{\gd^2})w^\gb + f(w^\gb)<0 \qtxt{in }\gw^\gb,$$
because in $\gw^\gb$: $\gd(\gx)=\dist(\gx,\bdw)<\dist((\xi_1+\gb,\xi'), \bdw)$. In addition  $w^\gb=0$ on $\prt \gw^\gb\cap \Gw$ and $w^\gb$ is bounded on $\prt \gw^\gb\cap \bdw$ while $u_1\to \infty$ as $\xi\to \prt \gw^\gb\cap \bdw$. Thus $w^\gb<u_1$ on $\prt\gw^\gb$.
By comparison principle $w^\gb< u_1$ in $\gw^\gb$. Letting $\gb\to 0$ we obtain $w\leq u_1$ in $\gw$. Therefore,
by (4*),
$$\limsup_{\xi\to\Gg_{1,a}} u_2/u_1\leq 1, \forevery a\in (0,1/2).$$
But, by assumption, $u_1\leq u_2$. Hence
\begin{equation}\label{u_2/u_1}
u_2/u_1\to 1  \qtxt{as }\; \xi\to \Gg_{1,a} \forevery a\in (0,1/2).
\end{equation}
This holds for every point $y\in \bdw$,  with rate of convergence independent of $y$. Therefore, for every $\ge>0$ there exists $s>0$ \sth
$$u_1\leq u_2 \leq (1+\ge)u_1 \qtxt{in }\Gw_s.$$
Our assumptions on $f$ imply that $(1+\ge)u_1$ is a supersolution of (NE). Therefore, by comparison principle, $u_1\leq u_2 \leq (1+\ge)u_1$ in $\Gw$. Letting $\ge\to 0$ we conclude that $u_1=u_2$.

\qed

\section{Proof of Theorem \ref{Main}} Throughout  this section we assume conditions (F1) and (F2) and the convexity of $f$. These conditions imply that $h$ is non-decreasing. As mentioned before (F2) implies (KO). We also assume that $\Gw$ is a bounded domain of class $C^2$ and that $\mu>0$. Other assumptions will be mentioned as needed.

Denote:
\begin{equation}\label{Gwgr}
\Gw_\gr:=\{x\in \Gw: \gd(x)<\gr\},
\end{equation}
\begin{equation}\label{tlphi}
\tl\phi(\gd):= h^{-1}(\gd^{-2}) ,\quad \gd>0.
\end{equation}
By  \eqref{con3} and \eqref{con2}, there exists $C>0 $ \sth
\begin{equation}\label{phi<tl}
h(\phi)\leq Ah(\tl \phi)\leq h(C \tl\phi) \qtxt{hence }\; \phi\leq  C\tl\phi.
\end{equation}

The proof of the theorem is based on several lemmas

\begin{lemma}\label{basic1}  Assume \eqref{con2}.
 	(i) There exists a maximal solution $U_{max}$ of (NE).
	 It satisfies
	 \begin{equation}\label{Umax}
	 h(U_{max}/a_1)\leq  a_0 \gd^{-2}= a_0  h(\tl\phi)
	 \end{equation}
$U_{max}$ is a large solution and
\begin{equation}\label{Uest1}
\rec{C}\phi<U^\Gw_f < U_{max}.
\end{equation}
\\
(ii) Every positive solution $u$ of (NE) satisfies
\begin{equation}\label{Uest2}
u \leq  \bar A \tl\phi
\end{equation}
where $\bar A$ is a constant depending only on the numbers $a_0, a_1$ in \eqref{Umax}.
\end{lemma}
\proof
In view of (F2) - which implies (F2') - the family of positive solutions of (NE) is uniformly bounded in every compact subset of the domain. Therefore there exists a maximal solution $U_{max}$ and, by (F2') and \eqref{tlphi}, it satisfies \eqref{Umax}.

$U^\Gw_f$ is a subsolution of (NE). The smallest solution between $U^\Gw_f$ and $U_{max}$ is a large solution of (NE). Therefore, by \eqref{u/phi}, we obtain  \eqref{Uest1}.

Let $u$ be an arbitrary positive solution of (NE). By \eqref{Umax} and \eqref{con2}  there exists a constant $ A'$  \sth
\begin{equation}\label{phi>U}
h(u/a_1)\leq a_0 \gd^{-2}= a_0 h(\tl\phi)\leq h( A'\tl\phi).
\end{equation}
Since $h$ is monotone this implies \eqref{Uest2} with $\bar A=A'a_1$.
\qed
\medskip

%
%

 Next we prove uniqueness of the large solution for every $\mu\geq 0$ \emph{under the conditions stated in Theorem \ref{Main}}. The proof is based on several lemmas.


\begin{lemma}\label{III.1}
(i) Let $D$ be a $C^2$ subdomain of $\Gw$ \sth $\Gg:=\prt D\cap \cap\bdw$ is the closure of a non-empty, relatively open subset of $\bdw$, say $O$.
Let $U$ be a positive  supersolution $U$ of (NE) in $D$. (Here $\gd(x)=\dist(x,\Gg)$.) If $U\to \infty $ as $x\to E$, for every $E\Subset O$ then,
\begin{equation}\label{mu<limsup}
 \mu\leq \limsup_{x\to E} h(U)\gd^2, \qtxt{for every compact $E\subseteq O$}.
\end{equation}
In particular, every large solution  satisfies the above inequality for $x\to \bdw$.

(ii) Using the notation of part (i): if $U$ is a positive supersolution of (NE) in $D$ then
\begin{equation}\label{mu<limsup'}
\mu-\rec{4}\leq \limsup_{x\to E} h(U)\gd^2,  \qtxt{for every compact $E\subseteq O$}.
\end{equation}
\end{lemma}


\proof (i) By negation, suppose that there exists $\ge$ positive \sth $\mu>   h(U)\gd^2 +\ge$ in  $D\cap \Gw_\rho$ for some $\gr>0$. Then, in this set,
$$0\leq -\Gd U -\mu\gd^{-2} U+ h(U) U< -\Gd U -\ge\gd^{-2}U.$$
Thus $U$ is $\Gd$-superharmonic in $D\cap \Gw_\rho$ and consequently has a (classical) measure boundary trace on $O$.
This contradicts the assumption.

(ii) If \eqref{mu<limsup'} is false, there exists $\ge>0$ and $\gr>0$ \sth
$$h(u)\gd^2< \mu-\rec{4}-\ge \qtxt{in }\; D\cap\Gw_\gr.$$
 Hence,
$$0\leq -\Gd u -\mu\gd^{-2}u+ h(u)u< -\Gd u - (\rec{4}+\ge)\gd^{-2}u \qtxt{in }\; D\cap\Gw_\rho.$$
It is known \cite[Thm. 5]{MMP} that if $\gg>\rec{4}$ there is no  local positive supersolution of the equation $0=-\Gd u -\gg\gd^{-2}u$. This contradicts the previous inequality

\qed

%

\begin{lemma}\label{III.2}
Suppose that $\Gw$ is radially symmetric: a ball, the exterior of a ball or an annulus. Then there exists $\bar b\in (0,1)$ \sth every radially symmetric large solution $U$ of (NE) satisfies
\begin{equation}\label{bf<U}
\bar b U_f^\Gw \leq  U \qtxt{in }\Gw.
\end{equation}


\end{lemma}

\proof First consider the case when $\Gw$ is a ball, say $B_R(0)$, or the exterior of a ball $B'_R(0)=\BBR^N\sms B_R(0)$.

By Lemma \ref{III.1} there exists a constant $c>0$ \sth
\begin{equation}\label{III.2.1}
\limsup_{\gd\to 0} h(U) /h(\tl\phi)>c
\end{equation}
for every radially symmetric large solution $U=U(\gd(x))$ where $\gd(x)=R-|x|$ in $B_R$ and $\gd(x)= |x|-R$ in $B'_R$.
By \eqref{con2.0},  there exists $b_0\in (0,1)$ and $\gd_0>0$, depending on $c$, \sth
\begin{equation*}
ch(\tl\phi) \geq h(b_0\tl\phi), \quad 0<\gd<\gd_0.
\end{equation*}
Therefore, by \eqref{III.2.1}, there exists a sequence $\{\gd_n\}$ converging to zero  \sth
\begin{equation}\label{III.2.1'}
h(U(\gd_n))\geq ch(\tl\phi(\gd_n)) \geq h(b_0\tl\phi(\gd_n)).  
\end{equation}
Hence, by  \eqref{u/phi} and \eqref{phi<tl}, there exists $\bar b>0$ \sth for sufficiently large $n$,
\begin{equation}\label{III.2.2}
U(\gd_n) \geq b_0\tl\phi(\gd_n)\geq \frac{b_0}{C}\phi(\gd_n)\geq \bar b\,U^\Gw_f(\gd_n).
\end{equation}
As $U^\Gw_f$ is a subsolution we obtain \eqref{bf<U}.

Now, let  $\Gw$ be an annulus, $R_1<|x|<R_2$. Then, by Lemma \ref{III.1}, there exist sequences $\{ \gd_n\}$ where $\gd_n=\gd(x_n)=|x_n|- R_1\to 0$ and   $\{\gd'_n\}$ where $\gd'_n=\gd(x'_n)=R_2- |x'_n|\to 0$ and a constant $c>0$ \sth \eqref{III.2.1'} and \eqref{III.2.2} hold for both sequences.
Thus 
$$U(\gd_n)\geq \bar b U_f^\Gw(\gd_n), \quad  U(\gd'_n)\geq \bar b U_f^\Gw(\gd'_n).$$
Hence, by the comparison principle,
$$U(x) \geq \bar b  U_f^\Gw(x),  \qtxt{when }\;   \gd_n \leq |x| \leq \gd'_n, $$
for all sufficiently large $n$. This imolies \eqref{bf<U}.  

\qed

\begin{lemma}\label{uniq1}
Suppose that $\Gw$ is radially symmetric and let $\CF_{rad}$ denote the family of
large, radially symmetric solutions. Then $\CF_{rad}$ contains a maximal element $U_{rad}$ and  a minimal element $u_{rad}$.
\end{lemma}

\proof The existence of a r.s. large solution follows from the fact that $U_f^\Gw$ is a r.s. large subsolution. The existence of a maximal element follows from the fact that, by (F2'), $\sup \CF_{rad}$ finite everywhere in $\Gw$ and is itself a solution (see e.g. [MVbook,p.79]).  Obviously  $\sup \CF_{rad}$ is r.s..
The existence of a minimal element follows from Lemma \ref{III.2}. Indeed if $U_1,U_2\in \CF_{rad}$ then $\min(U_1,U_2)\in \CF_{rad}$ is a supersolution larger than the subsolution $\bar b U_f^\Gw$. Therefore there exists $U_3\in \CF_{rad}$ \sth $U_3< \min(U_1,U_2)$. By a well known result, this fact implies that $\inf\CF_{rad}$ is a solution (see  [MVbook]) which obviously belongs to $\CF_{rad}$.
\qed

%
%
%
%
%
%
%
%
%

\begin{lemma}\label{uniq2}
Let $\Gw$ be a ball or the exterior of a ball. Then $\CF_{rad}$ is a sigleton: the unique \rs large solution.
\end{lemma}

\proof  By negation suppose that $U_{rad}\neq u_{rad}$. By \eqref{Uest2} 
$$U_{rad}\leq \bar A \tl\phi.$$
By  \eqref{III.2.2}, there exists $b_0>0$ and a sequence $\gd_n\downarrow 0$ \sth,
$$b_0\tl\phi(\gd_n)\leq u_{rad}(\gd_n).$$
Hence 
$$U_{rad}(\gd_n)\leq M u_{rad}(\gd_n) \qtxt{where }\;M= \frac{\bar A}{b_0}.$$
The function $M u_{rad}$  is a supersolution of (NE). Hence, by the comparison principle,
\begin{equation}\label{un2.1}
U_{rad}\leq M u_{rad} \qtxt{in }\Gw.
\end{equation}

Adapting a trick introduced in \cite{MV-sub} we show that this leads to a contradiction.

Put
$$w=u_{rad} - \rec{2M}(U_{rad}-u_{rad}).$$
By \eqref{un2.1},
$$\tl w:= \frac{M+1}{2M}u_{rad} < w< u_{rad} $$
The convexity of $f$ and the assumption $f(0)=0$ imply that $\tl w$ is a subsolution of (NE).  On the other hand, using these facts it is easily verified that $w$ is a supersolution of (NE). It follows that there exists a radially symmetric large solution
strictly smaller than $u_{rad}$, which brings us to a contradiction.

\qed

\noindent\textit{Completion of proof of the theorem.}\\ [2mm]
Choose $R>0$ \sth at every point $P\in \bdw$ there exists a ball $B_R(P')$ \sth $\bar B_R(P')\cap \bar \Gw=P$. Let $u_R$ be the unique \rs large solution of (NE) in
$B'_R(0)=\BBR^N\sms B_R(0)$. (Of course, here $\gd(x)=\dist (x,\prt B'_R(0))$.) Denote, $$u^P_{R}(x):=u_R(x-P') \forevery x\in B'_R(P'):=\BBR^N\sms B_R(P')$$ and $$u^P_{R,\ge}(x):= u_R(x-P'-\ge\mathbf{n}_P)$$.

Since $\mu>0$, $u^P_{R,\ge}$  is a subsolution of (NE) in $\Gw$  and it is dominated by $U$ on $\bdw$.
Therefore, $u^P_{R,\ge}(x)\leq U(x)$ in $\Gw$. Letting $\ge\to 0$, we obtaim
\begin{equation}\label{III.3.4}
u^P_R\leq U \qtxt{in }\Gw.
\end{equation}
The function $$U^*:=\sup_{P\in \bdw} u_R^P$$
 is a subsolution of (NE) in $\Gw$ and 
 $$U^*(x)\to\infty \qtxt{as }\; \Gw\ni x\to \bdw.$$ 
 By \eqref{III.3.4}, $U^*\leq U$ for every large solution $U$ of (NE) in $\Gw$. Therefore the smallest solution of (NE) dominating $U^*$ in $\Gw$ is the minimal large solution, denoted by $u_{min}$.

Let $\{\gd_n\}$ and $b_0$ be as in \eqref{III.2.2} with respect to $\BBR^N\sms B_R(0)$. Then
$$u_{min}(x)\geq U^*(x)\geq b_0\tl\phi(\gd_n) \forevery x\in \Gw: \gd(x)=\gd_n.$$
This inequality and \eqref{Uest2} imply that 
\begin{equation}\label{U<u}
U_{max}\leq M u_{min}, \quad M:= \frac{\bar A}{b_0}.
\end{equation}
Consequently, by the same argument as in the proof of  Lemma \ref{uniq2}, we conclude that $ U_{max}=  u_{min}$.

\qed

\end{document}